\documentclass{gtmon_a}
\pdfoutput=1

\usepackage{color,graphicx,pinlabel}

\proceedingstitle{The Zieschang Gedenkschrift}
\conferencestart{5 September 2007}
\conferenceend{8 September 2007}
\conferencename{Conference in honour of Heiner Zieschang}
\conferencelocation{Toulouse, France}

\editor{Michel Boileau}
\givenname{Michel}
\surname{Boileau}

\editor{Martin Scharlemann}
\givenname{Martin}
\surname{Scharlemann}

\editor{Richard Weidmann}
\givenname{Richard}
\surname{Weidmann}

\title{Noncoherence of some lattices in ${\rm Isom}(\mathbb{H}^n)$}

\author{Michael Kapovich}
\givenname{Michael}
\surname{Kapovich}
\address{Department of Mathematics\\
University of California, Davis\\\newline
1 Shields Ave\\ 
CA 95616\\ USA
}
\email{kapovich@math.ucdavis.edu}
\urladdr{}

\author{Leonid Potyagailo}
\givenname{Leonid}
\surname{Potyagailo}
\address{UFR de Math\'ematiques\\ Universit\'e de
Lille 1\\\newline
59655 Villeneuve d'Ascq cedex\\ France\vspace{3pt}\\\newline
Department of Mechanics and Mathematics\\
Lomonosov Moscow State University\\\newline
 Vorob'evy Gory\\Moscow 119992, GSP-2\\
Russia}
\email{potyag@math.univ-lille1.fr}
\urladdr{}

\author{Ernest Vinberg}
\givenname{Ernest}
\surname{Vinberg}
\email{vinberg@zebra.ru}
\urladdr{}

\dedicatory{To the memory of Heiner Zieschang}

\volumenumber{14}
\issuenumber{}
\publicationyear{2008}
\papernumber{015}
\startpage{335}
\endpage{351}

\doi{}
\MR{}
\Zbl{}

\arxivreference{math/0608415} 


\received{14 August 2006}
\revised{}
\accepted{18 November 2006}
\published{29 April 2008}
\publishedonline{29 April 2008}
\proposed{}
\seconded{}
\corresponding{}
\version{}

\makeatletter
\def\cnewtheorem#1[#2]#3{\newtheorem{#1}{#3}[section]
\expandafter\let\csname c@#1\endcsname\c@dfn}
\makeatother

\let\xysavmatrix\xymatrix
\def\xymatrix{\disablesubscriptcorrection\xysavmatrix}
\AtBeginDocument{\let\bar\wbar\let\tilde\wtilde\let\hat\what}
\def\SetFigFont#1#2#3#4#5{\small}

\theoremstyle{definition}
\newtheorem{dfn}{Definition}[section]
\cnewtheorem{rem}[dfn]{Remark}
\cnewtheorem{conj}[dfn]{Conjecture}
\cnewtheorem{ques}[dfn]{Question}
\cnewtheorem{probl}[dfn]{Problem}

\theoremstyle{plain}
\cnewtheorem{theorem}[dfn]{Theorem}

\newtheorem{thm}[dfn]{Theorem}
\cnewtheorem{lem}[dfn]{Lemma}

\newtheorem{thmA}{Theorem}

\makeautorefname{thmA}{Theorem}

\cnewtheorem{prop}[dfn]{Proposition}
\cnewtheorem{cor}[dfn]{Corollary}

\newcommand{\Isom}{\operatorname{Isom}}

\def\H{{\mathbb H}}

\def\N{{\mathbb N}}

\def\ga{\gamma}

\def\Ga{{\Gamma}}

\def\De{\Delta}
\def\Del{\Delta}

\def\si{\sigma}

\def\si{\sigma}

\def\dim{{\rm dim}}
\def\isom{{\rm Isom}(\H^n)}
\def\iso4{{\rm Isom}(\H^4)}

\def\h3{{{\mathbb H}^3}}
\def\h4{{\H^4}}

\begin{document}

\begin{asciiabstract} 
We prove noncoherence of certain families of lattices in the isometry
group of the hyperbolic n-space for n greater than 3. For instance,
every nonuniform arithmetic lattice in SO(n,1) is noncoherent,
provided that n is at least 6.
\end{asciiabstract}

\begin{htmlabstract}
We prove noncoherence of certain families of lattices in the isometry
group of the hyperbolic n&ndash;space for n greater than 3. For
instance, every nonuniform arithmetic lattice in SO(n,1) is
noncoherent, provided that n is at least 6.
\end{htmlabstract}

\begin{abstract} 
We prove noncoherence of certain families of lattices in the isometry
group of the hyperbolic $n$--space for $n$ greater than $3$. For
instance, every nonuniform arithmetic lattice in $SO(n, 1)$ is
noncoherent, provided that $n$ is at least $6$.
\end{abstract}

\maketitle

\section{Introduction}\label{intro}

The aim of this paper is to prove noncoherence of certain families
of lattices in the isometry group $\isom$ of the hyperbolic $n$--space $\H^n$ ($n>3$).
We recall that a group $G$ is called \textit{coherent} if every finitely generated subgroup of $G$
is finitely presented. It is well known that all lattices in $\Isom(\H^2)$ and $\Isom(\H^3)$
are coherent. Indeed, it is easy to prove that every finitely generated Fuchsian group is finitely
presented. The coherence of $3$--manifold groups  was proved by P\,Scott \cite{Sc1}. First examples
of geometrically finite noncoherent discrete subgroups of $\iso4$
were constructed by the first and second author \cite{KaP} and the second author \cite{P1,P2}. An example of
noncoherent uniform lattice in $\iso4$ was given by Bowditch and Mess~\cite{BM}.

In what follows we will identify $\H^n$ with a connected component of
the hyperboloid
$$
\{ x: f(x)=-1\}\subset \R^{n+1},
$$
where $f$ is a real quadratic form of signature $(n,1)$ in $n+1$
variables. Then the group $\isom$ is identified with the index $2$ subgroup $O'(f, \R)\subset O(f, \R)$
preserving~$\H^n$.

Let $f$ and $g$ be quadratic forms on finite-dimensional vector spaces $V$ and $W$ over $\Q$.
It is said that $f$ \textit{represents} $g$ if the vector space $V$ admits an
orthogonal decomposition (with respect to $f$)
$$
V= V'\oplus V''
$$
so that $f|V'$ is isometric to $g$. In other words, after a change
of coordinates, the form $f$ can be written as
$$
f(x_1,...,x_n)=g(x_1,...,x_k)+ h(x_{k+1},..., x_{n})
$$
where $n=\dim(V)$ and $k=\dim(W)$.
Whenever $f$ represents $g$, a finite index subgroup  of
$O(g, \Z)$ is naturally embedded into $O(f, \Z)$.

The main result of this paper is:

\begin{thmA}\label{Theorem A} For every $n\ge 4$ and every rational
quadratic form $f$ of signature $(n,1)$ which represents the form
$$
q_3=-x_0^2+x_1^2+x_2^2+x_3^2,
$$
the lattice $O(f, \Z)$ is noncoherent.
\end{thmA}

\begin{cor}
\label{C1} For every $n\ge 4$ there are infinitely many
commensurability classes of nonuniform noncoherent lattices in
$\isom$.
\end{cor}

 We refer the reader to \fullref{proofA} for the discussion of
uniform lattices. By combining \fullref{Theorem A} with some standard facts on rational quadratic forms,
we prove:

\begin{thmA}\label{Theorem B} For $n\ge 6$ every nonuniform arithmetic lattice in $\isom$ is noncoherent.
\end{thmA}

As a by-product of the proof, in \fullref{rational}, we obtain a
simple proof of the following result of independent interest
(which was proven by Agol, Long and Reid \cite{ALR} in the case $n=3$). Recall that a subgroup of
a group $\Gamma$ is called \textit{separable} if it can be represented
as the intersection of a family of finite index subgroups of
$\Gamma$. For instance, separability of the trivial subgroup is nothing else than
residual finiteness of $\Gamma$.

\begin{thmA}\label{Theorem C}
In every nonuniform arithmetic lattice in $\isom$ ($n\le 5$),
every geometrically finite subgroup is separable.
\end{thmA}

We refer the reader to Bowditch \cite{Bowditch} for the definition of
geometrically finite discrete subgroups of $\isom$. Recall only that
every discrete group which admits a finitely-sided convex
fundamental polyhedron is geometrically finite.

In \fullref{non-arithmetic} we adopt the method of Gromov and Piatetski-Shapiro \cite{GP} to
obtain examples of nonarithmetic noncoherent lattices:

\begin{thmA}\label{Theorem D}  For each $n\ge 4$ there exist both
uniform and nonuniform noncoherent nonarithmetic lattices in
$\isom$.
\end{thmA}

The above results provide a strong evidence for the negative answer
to the following question in the case of nonuniform lattices:

\begin{ques}[D\,Wise]\label{question}
 \label{DWise} Does there exist a coherent lattice in
$\isom$ for any $n>3$?
\end{ques}

In \fullref{spec} we provide some tentative evidence for the
negative answer to this question in the uniform case as well.

Our proof of the noncoherence in the nonuniform case is different from
the one by Bowditch and Mess \cite{BM}: The finitely generated infinitely
presented subgroup that we construct is generated by four subgroups
stabilizing 4 distinct hyperplanes  in $\H^n$, while in the
construction used in \cite{BM} two hyperplanes were enough. Direct
repetition of the arguments used in \cite{BM} does not seem to work in
the nonuniform case.

\medskip
{\bf Acknowledgements}\qua During this work the first author was
partially supported by the NSF grants DMS-04-05180 and DMS-05-54349.
A part of this work was done when the first and the second authors
were visiting the Max Planck Institute for Mathematics in Bonn.
The work of the third author was partially supported by the SFB 701
at Bielefeld University.

The second author is deeply grateful to Heiner Zieschang who was his
host during his Humboldt Fellowship at the University of Bochum in
1991--1992. The third author is thankful to Heiner for many years of
his generous friendship.

\section{Preliminaries}
\label{prelim}

We refer the reader to Kapovich \cite{Kapovich} and Maskit \cite{Maskit} for the basics of discrete groups
of isometries of the hyperbolic spaces $\H^n$.

\medskip {\bf Notation}\qua Given a convex polyhedron
$Q\subset \H^n$ let $G(Q)$ denote the subgroup of $\isom$ generated
by the reflections in the walls of $Q$.

We will frequently use the quadratic forms
$$
q_n=-x_0^2+x_1^2+...+x_n^2.
$$
Let $f$ be a quadratic form
$$
f=\sum_{i, j} a_{ij} x_i x_j
$$
defined over a number field $K\subset \R$, and  $\si$ be an embedding $K\to \R$.
Then $f^\si$ will denote the form
$$
\sum_{i, j} \si(a_{ij}) x_i x_j.
$$
%If $N$ is a normal subgroup in a finite index subgroup of a group
%$G$, then we say that $N$ is \textit{virtually normal} in $G$.

\subsection{Arithmetic groups}
\label{ar}

Let $f$ be a quadratic form of signature $(n,1)$ in $n+1$
variables with coefficients in a totally real algebraic number field
$K\subset \R$ satisfying the following condition:
\begin{equation}\label{*}\begin{array}{l}
\mbox{For every nontrivial (ie, different from the identity)
embedding $\si \co K\to \R$}\\ \mbox{the quadratic form $f^\si$ is
positive definite.}\end{array}\tag{$*$}
\end{equation}
Below we discuss discrete subgroups of $\isom$ defined using the form $f$.
Let $A$ denote the ring of integers of $K$.
We define the group $\Gamma:=O(f, A)$ consisting of matrices
with entries in $A$ preserving the form $f$. Then $\Gamma$ is
a discrete subgroup of $O(f, \R)$. Moreover, it is a \textit{lattice},
ie, its index 2 subgroup
$$
\Gamma'=O'(f, A):=O(f, A)\cap O'(f, \R)$$
 acts on $\H^n$ so that $\H^n/\Gamma'$ has finite volume.
Such groups $\Gamma$ (and subgroups of $\isom$
commensurable to them) are called \textit{ arithmetic
subgroups of the simplest type} in $O(n,1)$; see Vinberg and Shvartsman \cite{VS88}.

\begin{rem}
If $\Gamma\subset \isom$ is an arithmetic lattice so that
either $\Gamma$ is nonuniform or $n$ is even, then it follows from
the classification of rational structures on $\isom$ that $\Gamma$
is commensurable to an arithmetic lattice of the simplest type. For odd $n$ there
is another family of arithmetic lattices given as
the groups  of units of appropriate skew-Hermitian forms over
quaternionic algebras. Yet other families of arithmetic lattices exist for $n=3$ and $n=7$. See, for example,
Vinberg and Shvartsman \cite{VS88} or Millson and Li \cite{MillsonLi}.
\end{rem}

A lattice $\Gamma\subset \isom$ is called \textit{uniform} if
$\H^n/\Gamma$ is compact and \textit{nonuniform} otherwise. An
arithmetic lattice $O(f, A)$ of the simplest type is nonuniform if and only if $K=\Q$ and
$f$ is \textit{isotropic}, ie, there exists a nonzero vector
$v\in \Q^{n+1}$ such that $f(v)=0$.

Meyer's theorem (which follows from the Hasse--Minkowski principle;
see \cite[pp\,61--62]{BoSh} or \cite[Corollary 1, p\,75]{Ca})
states that every indefinite rational quadratic form of rank~$\ge 5$
is isotropic. Thus, for each rational quadratic form $f$ of
signature $(n,1)$, $n\ge 4$, the lattice $O(f, \Z)$ is nonuniform.
Conversely, every nonuniform arithmetic lattice in $\isom$ is
commensurable to $O(f, \Z)$, where $f$ is a rational quadratic form.

In particular, the groups $O'(q_n, \Z)\subset \isom$ are nonuniform
arithmetic lattices. The group $O'(q_3, \Z)$ coincides with the
group $G(\Delta)$, where $\Delta\subset \H^3$ is the simplex with
the Coxeter diagram
\begin{center}
\labellist
\small\hair2pt
\pinlabel $1$ [b] at 3 6
\pinlabel $2$ [b] at 74 6
\pinlabel $3$ [b] at 146 6
\pinlabel $4$ [b] at 219 6
\endlabellist
\includegraphics[width=2in]{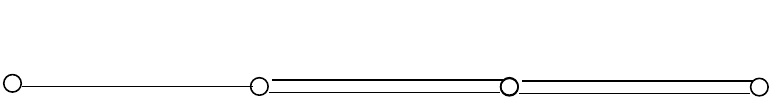} 
\label{f1.fig}
\end{center}
 (see \cite[Chapter 6, 2.1]{VS88} and references therein).

\begin{lem}
\label{fiber} The group $G(\Del)$ contains a finite index subgroup
$\Gamma$ such that $\H^3/\Ga$ fibers over the circle.
\end{lem}

\begin{proof} Let $v_{4}\in \De$ denote the (finite) vertex of $\Del$
disjoint from the $4$--th  face. Consider the union $O$ of the images
of $\Del$ under the stabilizer of $v_{4}$ in $G(\Del)$. Then $O$ is
a regular right-angled ideal hyperbolic octahedron in $\H^3$
\cite{Po-Vi,VS88}. The group $G(\Delta)$ contains $G(O)$ as a
finite index subgroup. It is well known that $G(O)$ is commensurable
with the fundamental group of the Borromean rings complement which
fibers over the circle \cite{Th}. The property of being the
fundamental group of a surface bundle over the circle is hereditary
with respect to subgroups of finite index. Thus $G(\Del)$ contains a
subgroup $\Gamma$ of finite index so that $\H^3/\Ga$ fibers over the
circle. \end{proof}

\subsection{Rational quadratic forms}
\label{rational}

The following proposition is well-known in the theory of rational
quadratic forms; see Cassels \cite[Exercise 8, Page 101]{Ca}. We present a proof
for the sake of completeness.

\begin{prop}
\label{forms} Let $f$ and $g$ be nonsingular rational quadratic
forms having respectively the  signatures $(r,s)$ and $(p,q)$ such
that $r\geq p$ and $s\geq q.$ If ${\rm rank}(f)-{\rm rank} (g)
\geq 3$ then $f$ represents $g$.

\end{prop}
\begin{proof}Recall that a rational quadratic form $f$ on a rational vector space
$V$ \textit{represents} $b\in \Q$ if there exists a vector
$v\in V\setminus \{0\}$ such that $f(v)=b$.  We use the following lemma.

\begin{lem}
\label{number} Suppose that $f$ is a nonsingular rational quadratic form in
$n\ge 4$ variables and $b$ is an
arbitrary nonzero rational number. 
\begin{enumerate}
\item[\rm(a)] If $f$ is positive definite and $b>0$ then $f$ represents $b$.

\item[\rm(b)] If $f$ is indefinite then $f$ represents $b$.
\end{enumerate}
\end{lem}
\begin{proof} The form
$$
F(y_1,...,y_n, y_{n+1}):=f(y_1,...,y_n)-b y_{n+1}^2
$$
is an indefinite nonsingular form of rank $\ge 5$.
By Meyer's theorem the form $F$ represents $0$. Hence
by \cite[Theorem 6, p\,393]{BoSh}, the form $f$ represents $b$.
 \end{proof}

 Let $n:=r+s,
k:=p+q$ be the ranks of $f$ and $g$ respectively. After changing
coordinates in $\Q^k$ we may assume that $g$ has the diagonal form
$$
g=b_1x_1^2+...+b_kx_k^2,
$$
where $b_i\in \Q_+$, if $i\le p$ and $b_i\in
\Q_-$, if $i>p$.

The form $f$ is isomorphic to
$b_1y_1^2+f_1(y_2,..., y_n)$ since $f$ represents $b_1$ by \fullref{number}.
By applying the same procedure to $f_1$ and arguing inductively we
obtain, after $k$ steps,
$$
f=b_1y_1^2+...+b_ky_k^2+f_k
$$
where $f_k$ is a form in $n-k$ variables. Note that the argument works
as long as $n-k\geq 3$. Indeed, if $n=k+3$ we will have
$$
f=b_1y_1^2+...+b_{k-1}y_{k-1}^2+f_{k-1}(y_k, y_{k+1}, y_{k+2},
y_{k+3})$$ and therefore we can apply the above argument the last
time to $f_{k-1}$. \end{proof}

We now use the above proposition to prove \fullref{Theorem C} stated
in \fullref{intro}.

\begin{proof} We will use the following result proven by P\,Scott in
\cite{Sc} for the convex--cocompact subgroups and by Agol, Long and Reid \cite{ALR} for the
geometrically finite subgroups:

\medskip{\sl Suppose that $P\subset \H^n$ is a right-angled polyhedron of
finite volume. Then every geometrically finite subgroup of $G(P)$ is
separable.}

\medskip Let $\Ga$ be a nonuniform arithmetic lattice in $\Isom(\H^k)$,
$k\le 5$. Then $\Ga$ is commensurable to $O(g, \Z)$ where $g$ is a
nonsingular rational quadratic form of signature $(k, 1)$.

According to \cite{Po-Vi} there exists a right-angled noncompact
convex polyhedron of finite volume $P^8\subset \H^8$. Moreover, the
group $G(P^8)$ is a finite index subgroup in $O'(q_8, \Z)$, see
\cite{VS88}. Since ${\rm rank}(q_8)-{\rm rank}(g)\ge 3$, it follows
that $q_8$ represents $g$, see \fullref{forms}. Hence we
have a natural embedding of a finite index subgroup of $\Gamma$ into
$G(P^8)$. As $P^8$ is right-angled, every geometrically finite
subgroup of $G(P^8)$ is separable. Since subgroup separability is
hereditary with respect to passing to a subgroup, we conclude that
every geometrically finite subgroup of $\Ga$ is separable. \end{proof}

\subsection{Hyperplane separability}
\label{separability}

In \fullref{non-arithmetic} we will need the following variation
on subgroup separability. Suppose that $\Gamma=O'(f, \Z)$ is an
arithmetic subgroup of $\isom$, where $f$ is a rational quadratic
form of signature $(n,1)$. Let $V_i\subset \R^{n+1}$, $i=0, 1,...,k$
be rational vector subspaces of codimension 1, so that
$V_i\otimes \R$ intersects
$\H^n$ along the hyperplane $H_i$, $i=0, 1,...,k$. We assume that
\begin{equation}\label{eq1}
H_0\cap H_i=\emptyset, \quad i=1,...,k.
\end{equation}
The following proposition is a generalization of Long \cite{Long}; its proof follows the lines of the
proof of Margulis and Vinberg \cite[Lemma 10]{MV}.

\begin{prop}
\label{separ}
There exists a finite index subgroup $\Gamma'\subset \Gamma$ so that for every $\ga\in \Ga'$
either $\ga(H_0)=H_0$ or
$$
\ga(H_0)\cap (H_0\cup H_1\cup ...\cup H_k)=\emptyset .
$$
\end{prop}
\begin{proof} Let $(\cdot , \cdot)$ denote the symmetric bilinear form on $\R^{n+1}$ corresponding to $f$.
Suppose that $V, V'\subset \R^{n+1}$ are codimension 1 vector subspaces which intersect
$\H^n$ along hyperplanes $H, H'$. Let $e, e'\in \R^{n+1}$ be nonzero vectors orthogonal to
$V, V'$ respectively. Then $H$ intersects $H'$ transversally iff
$$
|(e, e')|< \sqrt{( e, e) (e', e')}.
$$
For each $V_i$ ($i=0, 1, ...,k$) choose an orthogonal primitive integer vector $e_i$.
Then \eqref{eq1} implies that
$$
|(e_i, e_0)|\ge \sqrt{( e_i, e_i) (e_0, e_0)}, \quad i=1,...,k.
$$
Choose a natural number  $N$   which is greater than
$$
2\mskip-3mu\max_{i=0, 1, ...,k}  |( e_0, e_i)|.
$$
Let $\Gamma'=\Gamma(N)$ denote the level $N$ congruence subgroup in $\Gamma$, ie,
the kernel of the natural homomorphism
$$
\Gamma \to GL(n+1, \Z/N\Z).
$$
Then for every $\gamma\in \Gamma'$, $i=0, 1, ...,k$,
$$
(\gamma(e_i), e_0)\equiv (e_i, e_0) \hbox{~~(mod}~~ N)
$$
and therefore either
\begin{align*}
|(\ga(e_i), e_0)|&=|(e_i, e_0)|
\\
|(\gamma(e_i), e_0)| &> |(e_i, e_0)| \ge \sqrt{( e_i, e_i) (e_0, e_0)}=
\sqrt{( \gamma(e_i), \gamma(e_i)) (e_0, e_0)},\tag*{\hbox{or}}
\end{align*}
hence either $\gamma(H_0)=H_i$ or $\gamma(H_0)\cap H_i=\emptyset$.

Lastly, we have to ensure that $\ga(e_0)\ne \pm e_i$ for $i=1,...,k$ and all $\ga\in \Ga'$.
This is achieved by taking $N$ which does not divide some nonzero entries of
$e_0+e_i$ and of $e_0- e_i$ for all $i=1,...,k$.
\end{proof}

\subsection{A construction of noncoherent groups}
\label{nonco}

Let $L\subset \Isom(\H^3)$ be a subgroup commensurable to the reflection group $G(\Delta)$
defined in \fullref{ar}. We embed $\H^3$ in $\H^4$ as a
hyperplane $H$ and  naturally extend the action of $L$ from $H$ to
$\H^4$. Let $p_1, p_2\in \partial H$ be distinct parabolic points of
$L$. Let $\Pi_1, \Pi_2$ be perpendicular
 hyperplanes in $\H^4$ which are parallel to $H$ and asymptotic to $p_1, p_2$,
 respectively. Let $\tau_i$ denote the (commuting) reflections in
 $\Pi_i, i=1, 2$. Set $\tau_3:=\tau_1\tau_2$.
 Let $G$ denote the subgroup of $\iso4$ generated by $L, \tau_1, \tau_2$.

\begin{thm}
{\rm \cite{KaP}}\qua \label{kptheorem}
For every choice of the group  $L$, hyperplane $H$, points $p_1, p_2$ and hyperplanes
$\Pi_1, \Pi_2$ as above, the group $G$ is noncoherent.
\end{thm}

We will need the following:

\begin{cor}
\label{kpcor}
Suppose that $L_0, L_1, L_2, L_3$ are arbitrary finite index subgroups in
$$
L,~~ \tau_1 L\tau_1, ~~\tau_2 L \tau_2, ~~\tau_3 L \tau_3, \quad \hbox{respectively}.$$
Then the subgroup $S$ of $G$ generated by $L_0, L_1, L_2, L_3$
is noncoherent.
\end{cor}
\begin{proof} The intersection
$$
L':=L_0\cap \tau_1 L_1\tau_1\cap \tau_2 L_2 \tau_2\cap \tau_3 L_3 \tau_3$$
is a finite index subgroup in $L$. Let $S'$ denote the subgroup of $S$ generated by
\begin{equation}\label{eq2}
L', ~~\tau_1 L'\tau_1, ~~\tau_2 L' \tau_2, ~~\tau_3 L' \tau_3.
\end{equation}
It is clear that $S'$ has index $4$ in the group generated by $L', \tau_1, \tau_2$.
Since the latter is noncoherent by \fullref{kptheorem}, it follows that $S'$,
and thus $S$, is noncoherent as well. \end{proof}

\begin{figure}[ht!]
\begin{center}
\begin{picture}(0,0)%
\includegraphics[scale=.8]{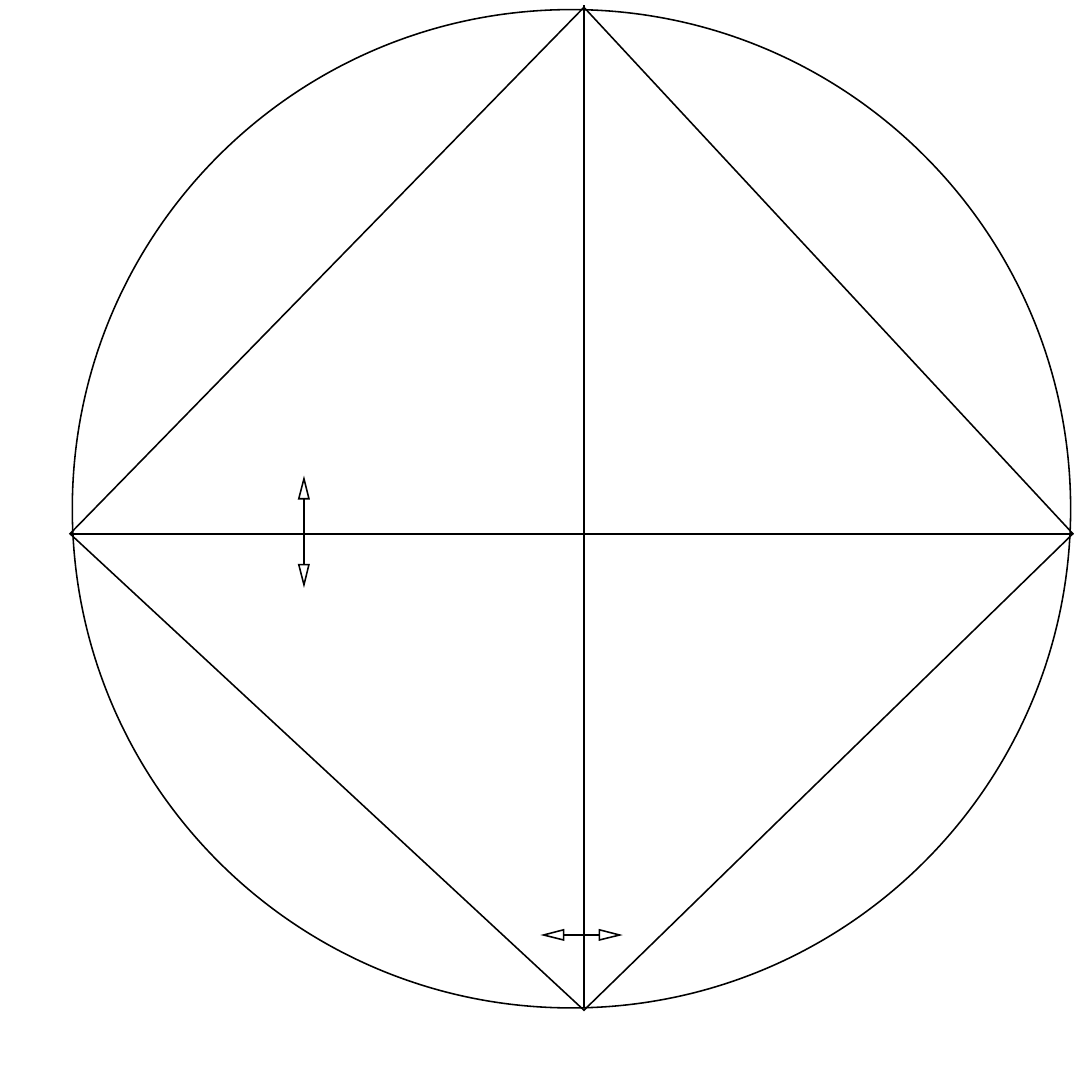}%
\end{picture}%
\setlength{\unitlength}{3157.6sp}%
\begingroup\makeatletter\ifx\SetFigFont\undefined%
\gdef\SetFigFont#1#2#3#4#5{%
  \reset@font\fontsize{#1}{#2pt}%
  \fontfamily{#3}\fontseries{#4}\fontshape{#5}%
  \selectfont}%
\fi\endgroup%
\begin{picture}(5147,5162)(3566,-5486)
\put(3566,-2894){\makebox(0,0)[lb]{\smash{\SetFigFont{10}{12.0}{\rmdefault}{\mddefault}{\updefault}{\color[rgb]{0,0,0}  $p_2$  }%
}}}
\put(6445,-3696){\makebox(0,0)[lb]{\smash{\SetFigFont{10}{12.0}{\rmdefault}{\mddefault}{\updefault}{\color[rgb]{0,0,0}  $\Pi_1$  }%
}}}
\put(5091,-3085){\makebox(0,0)[lb]{\smash{\SetFigFont{10}{12.0}{\rmdefault}{\mddefault}{\updefault}{\color[rgb]{0,0,0}  $\tau_2$  }%
}}}
\put(6274,-5440){\makebox(0,0)[lb]{\smash{\SetFigFont{10}{12.0}{\rmdefault}{\mddefault}{\updefault}{\color[rgb]{0,0,0}  $p_1$  }%
}}}
\put(5151,-3956){\makebox(0,0)[lb]{\smash{\SetFigFont{10}{12.0}{\rmdefault}{\mddefault}{\updefault}{\color[rgb]{0,0,0}  $H$  }%
}}}
\put(6835,-2723){\makebox(0,0)[lb]{\smash{\SetFigFont{10}{12.0}{\rmdefault}{\mddefault}{\updefault}{\color[rgb]{0,0,0}  $\Pi_2$  }%
}}}
\put(6425,-4698){\makebox(0,0)[lb]{\smash{\SetFigFont{10}{12.0}{\rmdefault}{\mddefault}{\updefault}{\color[rgb]{0,0,0}  $\tau_1$  }%
}}}
\put(5431,-1520){\makebox(0,0)[lb]{\smash{\SetFigFont{10}{12.0}{\rmdefault}{\mddefault}{\updefault}{\color[rgb]{0,0,0}$\tau_2(H)$}%
}}}
\put(6856,-1540){\makebox(0,0)[lb]{\smash{\SetFigFont{10}{12.0}{\rmdefault}{\mddefault}{\updefault}{\color[rgb]{0,0,0}$\tau_3(H)$}%
}}}
\put(5502,-2291){\makebox(0,0)[lb]{\smash{\SetFigFont{10}{12.0}{\rmdefault}{\mddefault}{\updefault}{\color[rgb]{0,0,0}  ${\mathbb H}^4$  }%
}}}
\put(7021,-3976){\makebox(0,0)[lb]{\smash{\SetFigFont{10}{12.0}{\rmdefault}{\mddefault}{\updefault}{\color[rgb]{0,0,0}  $\tau_1(H)$  }%
}}}
\end{picture}
\end{center}
\caption{}
\label{f2.fig}
\end{figure}

\begin{rem}
Note that the groups in \eqref{eq2} have the invariant hyperplanes
$H$, $\tau_1(H)$, $\tau_2(H)$, $\tau_3(H)$, respectively.
See \fullref{f2.fig}, where we use the projective model of $\H^4$.
\end{rem}

\section{Construction of noncoherent arithmetic lattices}
\label{proofA}

\begin{proof}[Proof of \fullref{Theorem A}] Our strategy is to embed a noncoherent
group $G$  (of the type described in \fullref{nonco}) into the
lattice $O(f, \Z)$. Then it would follow that $O(f, \Z)$ is
noncoherent.

Let $q_3$ be the quadratic form of rank 4 on the rational vector
space $U$ as in \fullref{prelim}. Then $O'(q_3, \Z)=G(\Delta)$,
see \fullref{ar}. We can change the coordinates in $U$ to $y_i$
($i=1,2,3,4$) so that $q_3$ takes the form:
$$
g=2y_1y_2+ y^2_3+y_4^2.
$$
Let $\{e_1, e_2, e_3, e_4\}$ be the corresponding basis of $U$. Note
that the group $O(g, \Z)$ is commensurable to $O(q_3, \Z)$.

Let $(U, g)\to (V,f)$ be a rational embedding. Pick a nonzero
vector $e_5\in V$ orthogonal to $U$. Then
$$
a:=f(e_5)>0.
$$
Define a $5$--dimensional vector space  $W$ spanned by the
vector $e_5$ and $U$. Let $h$ be the restriction of the form $f$ to
$W$; hence we have $(U,g)\subset (W,h)\subset (V,f)$. It therefore
suffices to embed some noncoherent group $G$ (as in \fullref{nonco})  into the group $O'(h, \Z)$.

We let $(\cdot , \cdot)$
denote the bilinear form on $W$ corresponding to $h$. The space $W$
splits as the orthogonal direct sum $U\oplus \Q e_5$. We will consider $\H^4$ canonically
embedded in $W\otimes \R$ and identify $\H^3$ with the hyperplane
$H:=U\otimes \R\cap \H^4\subset \H^4$.

After replacing $e_2$ with $ae_2$ we obtain $(e_1, e_2)=a$. Set
$$
u_1:= e_1+e_5, u_2:= -e_2+ e_5.
$$
Thus
$$
(u_1, u_1)= (u_2, u_2)=a,\ \ (u_1, u_2)=(u_1, e_1)=(u_2, e_2)=0,\ \
(u_1, e_2)=-(u_2,e_1)=a.
$$
Let $U_i\subset W$ ($i=1, 2$) be the $4$--dimensional vector subspace
orthogonal to $u_i$. Since $a>0$, it follows that each $U_i\otimes \R$ ($i=1,
2$) intersects $\H^4$ along a hyperplane
 $\Pi_i$. The reflection
$$
\tau_i\co w\mapsto w- 2\frac{(w,u_i)}{(u_i, u_i)} u_i
$$
in the subspace $U_i$ is represented by a matrix with integer
coefficients in the basis $\{e_1,...,e_5\}$. Thus $\tau_i\in O'(h,
\Z)$, $i=1, 2$.

Because $g(e_i)=0$, the vector $e_i$ corresponds to a parabolic
point $p_i\in \partial \H^4$ of the group $O(g, \Z)$, $i=1, 2$.
Since $(u_1, u_2)=0$, it follows that $\Pi_1$ is perpendicular to
$\Pi_2$. Moreover, since $e_i\in U_i$, we conclude that $\partial
\Pi_i$ contains $p_i$, $i=1, 2$. Since
$$
(u_i, e_5)=\sqrt{(u_i, u_i) (e_5, e_5)},
$$
the hyperplane $\Pi_i$ is parallel to $H$; see the proof of
\fullref{separ}.

Let $L$ be a finite index subgroup of $O'(g, \Z)$ contained in
$O'(h, \Z)$. The group $G$, generated by $L, \tau_1, \tau_2$, is
contained in $O'(h, \Z)$. \fullref{kptheorem} then implies that
the lattice $O'(h, \Z)$ is noncoherent. Hence $O(f, \Z)$ is
noncoherent as well. \fullref{Theorem A} follows. \end{proof}

\begin{proof}[Proof of \fullref{C1}]  For any number $a\in \N$
consider the quadratic form
$$
f_a(x_0, x_1,...,x_n)=q_3(x_0, x_1, x_2, x_3) + ax_4^2 +
x_5^2+...+x_n^2.
$$
Each $f_a$ defines a nonuniform arithmetic lattice $O'(f_a, \Z)\subset \isom$.
Moreover, for infinitely many appropriately chosen primes $a$
these lattices are not commensurable. Since each form $f_a$ represents $q_3$,
\fullref{C1}  follows from \fullref{Theorem A}. \end{proof}

\begin{theorem}
\label{bm} For each $n\ge 4$ there exist uniform noncoherent
arithmetic lattices in $\isom$. Moreover, for each $n\ge 5$ there
are infinitely many commensurability classes of such lattices.
\end{theorem}
\begin{proof} The assertion is a rather direct corollary of
the result of Bowditch and Mess \cite{BM}, but we present a proof
for the sake of completeness. We start with a review of the example
of Bowditch and Mess \cite{BM} which is a noncoherent uniform
arithmetic lattice in $\iso4$.

Consider the right-angled regular $120$--cell $D\subset \H^4$. It is a
compact regular polyhedron; see for instance Davis \cite{Davis} or Vinberg and Shvartsman \cite{VS88}. It
appears that it was first discovered by Schlegel in 1883
\cite{Schlegel}, who was interested in classifying \textit{honeycombs}
in the spaces of constant curvature; see Coxeter \cite{Coxeter}.

Each facet of $D$ is a right-angled regular dodecahedron. Let
$\Gamma=G(D)\subset \iso4$ be the reflection group determined by
$D$. The group $\Gamma$ is commensurable to $O(q, A)$, where $q(x_0,
x_1, x_2, x_3, x_4)$ is the quadratic form given by the matrix
\begin{equation}\label{eq3}
\left[
\begin{array}{ccccc}
1& -\cos(\pi/5) & 0 & 0 & 0 \\
-\cos(\pi/5) & 1& -1/2 & 0 & 0 \\
0& -1/2 & 1 & -1/2 & 0 \\
0& 0 & -1/2 & 1 & -\cos(\pi/5)\\
0 & 0 & 0& -\cos(\pi/5)&1
\end{array}
\right]
\end{equation}
and $A$ is the ring of integers of the field $K=\Q(\sqrt{5})$. Thus
$\Gamma$ is a (uniform) arithmetic lattice. Consider the facets
$F_1, F_2$ of $D$ which share a common $2$--dimensional face $F$. There
is a canonical isomorphism $\varphi\co G(F_1)\to G(F_2)$ fixing $G(F)$
elementwise. The reflection group $G(F_1)$ contains a finite index
subgroup isomorphic to the fundamental group of a hyperbolic
$3$--manifold $M^3$ which fibers over $S^1$; see Thurston \cite{Th}. Let
$N_1\subset \pi_1(M^3)$ be a normal surface subgroup and set
$N_2:=\varphi(N_1)\subset G(F_2)$. In particular, both $N_1, N_2$
are finitely generated. On the other hand, $N_i\cap G(F)$ is a free
group $E$ of infinite rank, $i=1, 2$. One then verifies that the
subgroup of $\Gamma$ generated by $N_1$ and $N_2$ is isomorphic to $
N_1*_E N_2$ and therefore is not finitely presentable \cite{Ne}.
Hence $\Gamma$ is a noncoherent uniform arithmetic lattice in
$\iso4$.

In order to construct lattices in $\isom$ consider the quadratic
forms
$$
f_a(x_0, x_1,..., x_n)=q(x_0, x_1, x_2, x_3, x_4)+ ax_5^2+
x_6^2+...+x_n^2,
$$
 where $a\in \N$ are primes. Since $q^\si$ is
positive definite for the (unique) nontrivial embedding $\sigma \co
K\to \R$, it follows that each $O(f_a, A)$ is a uniform arithmetic
lattice in $O(f_a, \R)$. As in the noncompact case, the groups
$O(f_a, \R)$ are not commensurable for infinitely many primes $a$.
As $O(q, A)\subset O(f_a, A)$, the assertion follows. \end{proof}

\begin{rem}
\label{R}
 Clearly, the subgroup generated by any finite index
subgroups of $G(F_1)$ and $G(F_2)$ is noncoherent as well.
\end{rem}

\begin{rem}
The above construction produces only one commensurability class of
 noncoherent  lattices in $\iso4$.
Using noncommensurable arithmetic lattices in $\iso4$ containing
$G(F_1)$, one can construct infinitely many commensurability
clas\-ses of uniform noncoherent arithmetic lattices in $\iso4$.
\end{rem}

\begin{proof}[Proof of \fullref{Theorem B}] Let $\Gamma$ be  a nonuniform arithmetic lattice in $\isom$
where $n\geq 6$. Then $\Gamma$ is commensurable to $O(f, \Z)$ for
some rational form $f$ of signature $(n,1)$.   Since $n+1\ge 7$ and
$q_3$ has rank $4$, it follows from \fullref{forms} that $f$
represents $q_3$. Therefore, by \fullref{Theorem A}, the group $O(f,\Z)$ is
noncoherent. Thus $\Gamma$ is noncoherent as well. \end{proof}

\section{Nonarithmetic noncoherent lattices}
\label{non-arithmetic}

\begin{proof}[Proof of \fullref{Theorem D}]  We produce these noncoherent examples by using the
construction of nonarithmetic lattices in $\Isom(\H^n)$ due to
Gromov and Piatetski-Shapiro \cite{GP}. We begin with a review of
their construction.

Let $f$ be a quadratic form of signature $(n-1,1)$ in $n$ variables
with coefficients in a totally real algebraic number field $K\subset
\R$. Let $A$ denote the ring of integers of $K$. We assume that $f$
satisfies Condition \eqref{*} from \fullref{ar}.

We let $K_+$ denote the set of $a\in K$ such that for each embedding
$\si\co K\to \R$ we have $\si(a)>0$. For $a\in K_+$ we consider the
quadratic form
$$
h_a(x_0, x_1, ..., x_n)=f(x_0, x_1,..., x_{n-1})+ ax_n^2.
$$
It has signature $(n,1)$ and satisfies Condition (*). Then
$\Gamma_a:=O'(h_a, A)$ is a lattice in $\isom$. Similarly,
$\Gamma_0:=O'(f,A)$ is a lattice in $\Isom(\H^{n-1})$.

In what follows we will consider pairs of groups $\Gamma_a, \Gamma_1$, where $a\in \N$.
Observe that both groups contain the subgroup $\Gamma_0$.
Let $\Gamma_a'\subset \Gamma_a, \Gamma_1'\subset \Gamma_1$ be
torsion-free finite index subgroups such that
$$
\Gamma_1'\cap \Gamma_0= \Gamma_a'\cap \Gamma_0.
$$
We let $\Gamma_0'$ denote this intersection and set $M_1:=\H^n/\Gamma_1', M_a:=\H^n/\Gamma_a'$.

Without loss of generality (after passing to deeper finite index
subgroups), we may assume that $\H^{n-1}/\Gamma_0'$
isometrically embeds into $M_1$ and $M_a$ as a nonseparating totally geodesic hypersurface; see
Millson \cite{Millson}. Cut $M_1$ and $M_a$ open along these hypersurfaces.
The resulting manifolds $M_1^+,
M_a^+$ both have totally geodesic boundaries isometric to the
disjoint union of two copies of $M_0=\H^{n-1}/\Gamma_0'$.

Let $M$ be the connected hyperbolic manifold obtained by gluing
$M_1^+, M_a^+$ via the isometry of their boundaries. It is easy to
see that $M$ is complete. Then there exists a lattice $\Gamma\subset
\Isom(\H^n)$ such that $M=\H^n/\Gamma$. Note that $M$ is compact iff
both $M_1, M_a$ are. It is proven in \cite{GP} that

\medskip 
{\sl $\Gamma$ is not arithmetic if and only if  $a$ is not a
square in $K$.}

\medskip
 Note that there exist infinitely numbers $a$ which are not squares in $K$.
Indeed, it is well known that square roots of prime numbers are linearly
independent over $\Q$. Therefore only finitely many of them belong to $K$.

We now prove \fullref{Theorem D} by working with specific examples.

\medskip
{\bf (1) Compact case}\qua For $n\ge 5$ take $K=\Q(\sqrt{5})$ and
consider the quadratic form
$$
f=q(x_0, x_1, x_2, x_3, x_4)+ x_5^2+...+x_{n-1}^2,$$ where the form
$q$ is given by the matrix \eqref{eq3}.
%Let $a$ be any prime number different from $5$.
The quadratic form $q$ yields a uniform arithmetic lattice $O'(q,
A)$ in $\Isom(\H^4)$ which is commensurable to the reflection group
$G(D)$ defined in the proof of \fullref{bm}. The group $O'(q,
A)$ is noncoherent according to \fullref{bm}. On the other
hand, by applying the Gromov--Piatetski-Shapiro construction to
$f$ and taking any prime number $a\ne 5$, we
obtain a uniform nonarithmetic lattice in $\isom$ which contains
$O'(q, A)$ and, hence, is noncoherent.

It remains to analyze the case $n=4$. Take facets $F_1, F_2, F_3$ of
$D$ so that $F_1$ and $F_2$ intersect along a $2$--dimensional face and
$$
F_3\cap F_1=F_3\cap F_2=\emptyset.
$$
Then the group generated by the reflections in the facets of $F_1$
and $F_2$ is noncoherent; see the proof of \fullref{bm}.

By taking an appropriate finite index subgroup $\Gamma_1\subset
G(D)$, we obtain a hyperbolic $4$--manifold $M_1=\H^4/\Gamma_1$ which
contains embedded totally geodesic hypersurfaces $S_i$ corresponding
to the facets $F_i$, $i=1, 2,3$, so that
$$
S_3\cap S_1=S_3\cap S_2=\emptyset, \quad S_1\cap S_2\ne \emptyset .
$$
Now cut $M_1$ open along $S_3$ and apply the gluing construction of
Gromov and Piatetski--Shapiro. In this way one can obtain a
nonarithmetic compact hyperbolic manifold $M$ whose fundamental
group contains the subgroup of $\Gamma_1$ generated by some finite
index subgroups of $G(F_1)$ and $G(F_2)$ and, hence, is noncoherent
(see \fullref{R}).

\medskip
{\bf (2) Noncompact case}\qua For $n\ge 5$ take $K=\Q$ and consider the
quadratic form $f=q_{n-1}$. Taking any prime number for $a$, apply
the same argument as in the compact case.

Consider $n=4$. We will imitate the proof in the compact case.
However we will appeal to the results of \fullref{separability}
instead of using a particular fundamental domain.

Let $\Gamma:=O'(q_4, \Z)$. Clearly, $q_4$ represents  the form
$q_3$. Set $L:=O'(q_3,\Z)\subset \Gamma$, and let $\tau_1, \tau_2
\in O'(q_4,\Z)$ be the commuting reflections constructed in the
proof of \fullref{Theorem A}. Set $\tau_3:=\tau_1 \tau_2$ and $L_0:=L,
L_i:=\tau_i L \tau_i$, $i=1, 2, 3$.

By passing to any finite index subgroups $L_i'\subset L_i$, we obtain
a noncoherent subgroup $G'$ in $\Gamma$ generated by $\smash{L_i'}, i=0, 1,
2, 3$; see \fullref{kpcor}. Since $\Gamma$ is a linear group,
we can assume without loss of generality that $\Gamma$ is
torsion-free.
Let $H_0=H\subset \H^4$ be the $L$--invariant hyperplane.
Then $H_i:=\tau_i (H)$ ($i=1,2,3$) is the $L_i$--invariant hyperplane.

\begin{lem}
There exists a finite index subgroup $\Gamma'\subset \Gamma$ so that
for the groups $L_i':=L_i\cap \Gamma'$ we have:
\begin{enumerate}
\item $H_0/L_0'$ embeds as a hypersurface $S_0$ into $\H^4/\Gamma'$.

\item Let  $M^+$ denote the manifold obtained by cutting $\H^4/\Gamma'$
along $S_0$. Then $G'$ embeds into $\pi_1(M^+)$.
\end{enumerate}
\end{lem}
\begin{proof} We have to find a subgroup $\Ga'$ so that:

(a)\qua For all $\ga\in \Ga'$ either $\ga(H_0)=H_0$
or
$$
\ga(H_0) \cap( H_0\cup H_1 \cup H_2 \cup H_3)=\emptyset.
$$
(b)\qua For all $\ga\in \Ga'$,  the hyperplane $\ga(H_0)$ does not
separate the above hyperplanes from each other.

This is achieved by applying \fullref{separ} to
the hyperplanes $H_0, H_1, H_2, H_3$
and $H_4:=\Pi_1$, $H_5:=\Pi_2$.  \end{proof}

We now glue an appropriately chosen manifold $M_a^+$ along the
boundary of $M^+$. Let $M$ be the resulting complete hyperbolic
manifold. Then, as in the case $n\ge 5$, the fundamental group of
$M$ is nonarithmetic. On the other hand,
$$
G'\subset \pi_1(M^+)\subset \pi_1(M).
$$
Therefore $\pi_1(M)$ is noncoherent. \end{proof}

\section{Noncoherence and Thurston's conjecture}
\label{spec}

We recall the following conjecture:

\begin{conj}
[Thurston's virtual fibration conjecture] Suppose that $M$ is a hyperbolic
$3$--manifold of finite volume. Then there exists a finite cover over
$M$ which fibers over the circle.
\end{conj}

We expect that all lattices in $\isom$ are noncoherent for $n\ge
4$. Proving this for nonarithmetic lattices is clearly beyond our
reach.  Therefore we restrict to the arithmetic case. Even in this
case our discussion will be rather speculative. We restrict to the
arithmetic groups of the simplest type $\Gamma=O(f, A)$, where $f$
is a quadratic form on $V=K^{n+1}$ and $K\subset \R$ is a totally real
algebraic number field (see \fullref{ar}). Choose a basis
$\{e_0, e_1,...,e_n\}$ in which the form $f$ is diagonal:
$$
f=a_0 x_0^2+ a_1 x_1^2+... + a_n x_n^2.
$$
Here $a_0<0$, $a_1,...,a_n>0$ and for all nontrivial embeddings
$\si\co K\to \R$ we have $\si(a_i)>0$, for all $i=0, 1, ...,n$. To
simplify the discussion, we will assume that $\Gamma$ is uniform
(the nonuniform lattices were discussed in Theorems \ref{Theorem A} and \ref{Theorem B}).

For a $4$--element subset $I=\{0, i, j, k\}\subset \{0, 1,...,n\}$ let
$V_I\subset V$ denote the  linear span of the basis vectors $e_l, l\in I$. Set
$H_I:=V_I\otimes \R \cap \H^n$. Then $f|V_I$ determines a lattice $\Gamma_{I}$ in
$\Isom(H_I)$, which is naturally embedded into $\Ga$.  Assuming
Thurston's conjecture, up to taking finite index subgroups, each
$\Gamma_{I}$ contains an (infinite)  normal finitely generated surface subgroup
$N_{I}$. Moreover, by taking $I$ and $J$ such that $I\cap J$
consists of 3 elements, we obtain subgroups $\Gamma_{I}, \Gamma_{J}$
whose intersection is a Fuchsian group $F$. It now follows from the
 separability of  $F$ in $\Gamma$  (see Bowditch and Mess \cite{BM}, Long \cite{Long} or
 \fullref{separ} of this paper) that, after passing to certain finite index
subgroups $\Ga_I'\subset \Gamma_{I}, \Ga_J'\subset \Gamma_{J}$, we get
the inclusion
\begin{equation}\label{eq4}
\Ga_I' *_F \Ga_J'\subset \Gamma. 
\end{equation}
Set $N_{I}':= \Ga_{I}'\cap N_{I}, N'_{J}:= \Ga_J'\cap N_{J}$. Then
$E:=N_I'\cap N_J'$ is a free group of infinite rank. Now \eqref{eq4} implies
that $\Gamma$ is noncoherent since the subgroup
$$
N_I' *_E N_J'\subset \Gamma
$$
is finitely generated but not finitely presented \cite{Ne}.
Therefore we obtain:

\medskip{\sl Suppose that Thurston's conjecture holds for all compact
arithmetic $3$--manifolds. Then all uniform arithmetic lattices of the
simplest type in $\isom$, $n\ge 4$, are noncoherent.}

\medskip Thus we expect the negative answer to \fullref{question} asked by Dani Wise.

\bibliographystyle{gtart}
\bibliography{link}

\end{document}